%% file: 2000-16.tex
\documentclass{gtart}

\input gtoutput

\lognumber{122}
\volumenumber{4}\papernumber{16}\volumeyear{2000}
\pagenumbers{451}{456}
\proposed{Ronald Stern}
\seconded{Yasha Eliashberg, Ronald Fintushel}
\received{11 June 2000}
\accepted{21 November 2000}
\published{26 November 2000}
 
\usepackage{amsfonts}

\def\ZZ{\mathbb Z}

\newtheorem{thm}{Theorem}
\theoremstyle{definition}
\newtheorem*{defn}{Definition}

\newtheorem*{xpl}{Example}
\newtheorem*{ack}{Acknowledgment}

\begin{document}

\title{Diffeomorphisms, symplectic forms\\and Kodaira fibrations}
\author{Claude LeBrun}
\address{Department of Mathematics, SUNY at Stony Brook\\Stony Brook, 
NY 11794-3651, USA}
\email{claude@math.sunysb.edu}

\begin{abstract} 
As was recently pointed out by McMullen and Taubes \cite{mct}, there
are $4$--manifolds for which the diffeomorphism group does not act
transitively on the deformation classes of orientation-compatible
symplectic structures. This note points out some other $4$--manifolds
with this property which arise as the orientation-reversed versions of
certain complex surfaces constructed by Kodaira \cite{kodf}.  While
this construction is arguably simpler than that of McMullen and
Taubes, its simplicity comes at a price: the examples exhibited herein
all have large fundamental groups.
\end{abstract}
 
\asciiabstract{As was recently pointed out by McMullen and Taubes
[Math. Res. Lett. 6 (1999) 681-696], there are $4$--manifolds for
which the diffeomorphism group does not act transitively on the
deformation classes of orientation-compatible symplectic
structures. This note points out some other $4$--manifolds with this
property which arise as the orientation-reversed versions of certain
complex surfaces constructed by Kodaira [J. Analyse Math. 19 (1967)
207-215].  While this construction is arguably simpler than that of
McMullen and Taubes, its simplicity comes at a price: the examples
exhibited herein all have large fundamental groups.}


\asciikeywords{Symplectic manifold, complex surface, Seiberg-Witten invariants}

\maketitlepage

Let $M$ be a smooth, compact oriented $4$--manifold. 
If $M$ admits an orientation-compatible symplectic form, meaning a closed 
$2$--form $\omega$ such that $\omega \wedge \omega$ is an
orientation-compatible volume form, 
one might well ask whether the space of such forms is
connected. In fact, it is not difficult to construct examples
where the answer is negative. A more subtle question, however,
is whether the  group of orientation-preserving diffeomorphisms
$M\to M$ acts transitively on the set of connected components of
the orientation-compatible  symplectic structures of $M$.
As was recently pointed out by McMullen and
Taubes \cite{mct}, there are $4$--manifolds $M$ for which this
subtler question also has a negative 
answer. The purpose of the present note is 
to point out that  many examples of this interesting phenomenon
 arise from  certain  complex surfaces 
 with Kodaira fibrations.

A {\em Kodaira fibration} is by definition a holomorphic
submersion $f\co  M\to B$  from a compact complex surface to a 
compact complex curve, with  base $B$ and fiber 
$F_{z}= f^{-1}(z)$ both of genus $\geq 2$. (In $C^{\infty}$ terms, $f$ 
is thus 
 a locally trivial fiber
bundle, but nearby fibers of $f$ may well be 
non-isomorphic  as complex curves.)
 One says that $M$ is a 
{\em Kodaira-fibered surface}  if it admits
such a fibration $f$. Now any Kodaira-fibered surface 
 $M$ is algebraic, since $K_{M}\otimes f^{*}K_{B}^{\otimes \ell}$
 is obviously positive  
  for sufficiently large $\ell$.
  On the other hand,
  recall that a holomorphic map from 
a curve of lower genus to a curve of higher genus must be 
constant.\footnote{Indeed, by Poincar\'e duality, a continuous map
$h\co X\to Y$ of non-zero degree between compact oriented manifolds of
the same dimension
must induce inclusions $h^{*}\co H^{j}(Y,{\mathbb R}) \hookrightarrow  
H^{j}(X,
{\mathbb R})$ for 
all $j$. Such a map $h$ therefore  cannot exist whenever   
$b_{j}(X) < b_{j}(Y)$
 for some $j$.} If  $f\co  M\to B$ is a Kodaira fibration, 
 it follows that $M$ cannot contain any rational or elliptic
 curves, since composing $f$  with the inclusion  would
 result in a constant map, and the curve would therefore
 be contained in a fiber of $f$; contradiction. 
  The Kodaira--Enriques 
 classification 
 \cite{bpv} therefore tells us that $M$ is a minimal surface of 
 general type.  In particular,  the only
 non-trivial Seiberg--Witten invariants of the 
underlying oriented $4$--manifold $M$ are   \cite{morgan}
those  associated with the  
canonical and anti-canonical classes of $M$. Any 
orientation-preserving self-diffeomorphism of $M$ 
must therefore preserve $\{ \pm c_{1}(M)\}$. 

We have just seen that $M$ is of  K\"ahler type,  
so let  $\psi$ denote some K\"ahler form  on $M$, and observe that 
  $\psi$ is then of course   
a symplectic form compatible with 
the usual `complex' orientation of $M$. 
Let $\varphi$ be any area form on $B$, compatible with 
{\em its} complex orientation, and, for sufficiently small
$\varepsilon > 0$,  consider the closed $2$--form  
$$\omega = \varepsilon \psi- f^{*} \varphi.$$
Then 
 $$\frac{\omega \wedge \omega}{\varepsilon} = 
 -2(f^{*}\varphi )\wedge \psi +  \varepsilon\psi \wedge \psi =
  \left( \varepsilon -\langle 
f^{*}\varphi , \psi \rangle  \right) \psi \wedge \psi ,$$
where the inner product is taken with respect to the 
 K\"ahler metric corresponding to $\psi$. 
Now $ \langle 
f^{*}\varphi , \psi \rangle$ is a positive function, and,
because $M$ is compact, therefore has a  positive minimum.
Thus, 
for a  sufficiently small $\varepsilon > 0$, 
$\omega \wedge \omega$ is a   
volume form compatible with   the
{\em non-standard} orientation of $M$; or, in other words,  $\omega$
is  a symplectic form for  
the reverse-oriented $4$--manifold 
$\overline{M}$. For related constructions of 
symplectic structures on fiber-bundles, cf  \cite{mcsal}.

If follows that $\overline{M}$ 
carries a unique deformation class of almost-complex
structures compatible with $\omega$. One such almost-complex
structure can be constructed by considering 
the (non-holomorphic) orthogonal decomposition
$$TM= \ker (f_{*}) \oplus f^{*}(TB)$$
induced by the given K\"ahler metric, and then 
reversing the sign of the complex structure on the 
`horizontal' bundle  $f^{*}(TB)$.
The
first Chern class of the resulting  almost-complex structure
 is thus given by 
$$c_{1}(\overline{M}, \omega ) = c_{1}(M) - 4 (1-{\mathbf g}) F,
$$
where ${\mathbf g}$ is the genus of $B$, and where 
$F$ now denotes the Poincar\'e dual of a  fiber of  $f$.
For further discussion, cf \cite{kot2,leung,jpik}.

Of course,  the product $B\times F$ of two complex curves of
genus $\geq 2$ is certainly  Kodaira fibered, but such  
a product also admits orientation-reversing diffeomorphisms, and 
so, in particular,  has  signature
$\tau=0$. However,
as was first observed by 
Kodaira \cite{kodf}, one can  construct examples 
with $\tau > 0$ by taking {\em branched covers} of 
products; cf  \cite{atkodf,bpv}. 

\begin{xpl}
Let $C$ be a compact complex curve of genus $k\geq 2$, 
and let $B_1$ be a curve of genus ${\bf g}_1=2k-1$, obtained as an
unbranched double cover of $C$. Let $\iota \co  B_1\to B_1$ be the
associated non-trivial deck transformation, which is  a
free holomorphic 
involution of $B_1$. 
Let $p\co  B_2\to B_1$ be the 
unique unbranched cover of order $2^{4k-2}$ 
 with $p_{*}[\pi_{1}(B_2)]= \ker [\pi_{1}(B_1) \to H_{1}(B_1, \ZZ_{2})]$;
 thus $B_2$ is a complex curve of genus ${\bf g}_2=2^{4k-1}(k-1)+1$. Let 
 $\Sigma \subset B_2\times B_1$ be the union of 
 the graphs of $p$ and $\iota\circ p$. Then
 the homology class of $\Sigma$ is divisible by $2$. We may
 therefore construct a ramified double cover $M\to B_2\times B_1$
 branched over $\Sigma$. The projection $f_1\co  M\to B_1$
 is then a Kodaira fibration, with fiber $F_1$ of genus 
 $2^{4k-2}(4k-3)+1$.
 The projection $f_2 \co  
 M\to B_2$ is also a Kodaira fibration, 
 with fiber $F_2$ of genus $4k-2$. The signature  of this 
 doubly Kodaira-fibered complex surface is 
  $\tau (M) =2^{4k}(k-1)$.
  \end{xpl}
  
 We now axiomatize those properties of these examples
 which we will need. 
 
  \begin{defn}
  Let $M$ be a complex surface equipped with 
  two Kodaira fibrations
  $f_{j}\co M\to B_{j}$, $j=1,2$. 
  Let ${\mathbf g}_{j}$ denote the genus of $B_{j}$, and 
  suppose that the induced map 
  $$f_{1}\times f_{2}\co  M\to B_{1}\times B_{2}$$
  has degree $r > 0$. We will then say that 
  $(f_{1}, f_{2})$ is a  {\em Kodaira double-fibration}
  of $M$  if $\tau (M)\neq 0$ and 
  $$ ({\mathbf g}_2-1)\not | ~ r({\mathbf g}_1-1).$$
  In this case,  
  $(M, f_{1}, f_{2})$ will be called a {\em Kodaira doubly-fibered}
  surface.
\end{defn}

 Of course,  the  last hypothesis depends on the
 ordering of $(f_{1}, f_{2})$, and  is 
  automatically satisfied, for fixed $r$,  
  if  ${\mathbf g}_{2}\gg {\mathbf g}_{1}$.
  The latter   may
 always be arranged by simply 
 replacing $M$ and $B_{2}$ with suitable covering spaces. 
 
 Note that $r=2$ in the explicit examples given above.

 Given a  Kodaira doubly-fibered surface $(M, f_{1}, f_{2})$, 
 let $\overline{M}$ denote $M$ equipped with the non-standard orientation, and 
 observe that we now have two different symplectic 
 structures on $\overline{M}$ given by 
 \begin{eqnarray*}
 	\omega _{1}  & = &   \varepsilon \psi- f^{*}_{1} \varphi_{1}\\
 	\omega _{2}  & = &   \varepsilon \psi- f^{*}_{2} \varphi_{2}
 \end{eqnarray*}
for any given area forms $\varphi_{j}$ on $B_{j}$ and 
any sufficiently small $\varepsilon >0$.

 \begin{thm}\label{one}
 Let $(M, f_{1}, f_{2})$  be any 
 Kodaira doubly-fibered complex surface. Then for 
  any self-diffeomorphism $\Phi\co  M\to M$, the 
 symplectic structures $\omega_{1}$ and $\pm \Phi^{*}\omega_{2}$ are
 deformation inequivalent.
 \end{thm}
 
 That is, $\omega_{1}$, $-\omega_{1}$,$\Phi^{*}\omega_{2}$, and
  $-\Phi^{*}\omega_{2}$ are always in 
 different path components of the closed, non-degenerate 
 $2$--forms on $\overline{M}$. (The fact that 
 $\omega_{1}$ and  $-\omega_{1}$ are deformation inequivalent
 is due to a general result of Taubes \cite{taubes2}, and 
 holds for any symplectic $4$--manifold with 
   $b^{+}> 1$ and $c_{1}\neq 0$.)

 Theorem \ref{one} is actually a corollary of 
  the following result:

 \begin{thm}
 Let $(M, f_{1}, f_{2})$  be any 
 Kodaira doubly-fibered complex surface. Then for 
  any self-diffeomorphism $\Phi\co  M\to M$,
  $$\Phi^{*}[c_{1}(\overline{M}, \omega_{2})]\neq \pm 
  c_{1}(\overline{M}, \omega_{1}).$$
 \end{thm}
 
 \begin{proof}
 Because $\tau (M)\neq 0$, any self-diffeomorphism of 
 $M$ preserves orientation. Now $M$ is a minimal 
 complex surface of general type, and hence, 
 for the standard `complex' orientation of $M$, 
 the only Seiberg--Witten basic classes  \cite{morgan}  are 
 $\pm c_{1}(M)$. Thus any self-diffeomorphism $\Phi$ of 
 $M$ satisfies
 $$\Phi^{*}[c_{1}(M)]= \pm c_{1}(M).$$
 
 Letting $F_{j}$ be the Poincar\'e dual of the fiber of
 $f_{j}$, and letting ${\mathbf g}_{j}$ denote the 
 genus of $B_{j}$, we have
 $$c_{1}(\overline{M}, \omega_{j})= c_{1}(M)+ 4({\mathbf g}_{j}-1) F_{j}$$
 for $j=1,2$. 
 The adjunction formula therefore tells us that 
 $$[c_{1}(\overline{M}, \omega_{j})]\cdot [c_{1}(M)]=
 (2\chi + 3\tau )(M) -2\chi (M) = 3\tau (M) \neq 0,$$
 where the intersection form is computed with respect to
 the `complex' orientation of $M$. 
 
 If we had a diffeomorphism $\Phi \co  M\to M$ with 
 $\Phi^{*}[c_{1}(\overline{M}, \omega_{2})]= \pm 
  c_{1}(\overline{M}, \omega_{1})$, this computation would tell us that 
  that 
  $$\Phi^{*}[c_{1}(M)]= c_{1}(M) ~
  \Longrightarrow ~ \Phi^{*}[c_{1}(\overline{M}, \omega_{2})]=  
  c_{1}(\overline{M}, \omega_{1})$$
  and that 
  $$\Phi^{*}[c_{1}(M)]= -c_{1}(M) ~
  \Longrightarrow ~ \Phi^{*}[c_{1}(\overline{M}, \omega_{2})]=  
  -c_{1}(\overline{M}, \omega_{1}).$$
  In either case, we would then have 
  $$4({\mathbf g}_{1}-1) F_{1} =c_{1}(\overline{M}, \omega_{1})-c_{1}(M)=
  \pm \Phi^{*}[c_{1}(\overline{M}, \omega_{2})-c_{1}(M)]
  =\pm 4({\mathbf g}_{2}-1)\Phi^{*}( F_{2}).
  $$
  On the other hand, $F_{1}\cdot F_{2}=r$, so intersecting 
  the previous formula with $F_{2}$ yields  
  $$4({\mathbf g}_{1}-1)r= 4({\mathbf g}_{1}-1) F_{1}\cdot F_{2}= 
  4({\mathbf g}_{2}-1)
  [\pm \Phi^{*}( F_{2}) \cdot F_{2}],$$
  and hence   
  $$({\mathbf g}_{2}-1)~|~ r({\mathbf g}_{1}-1),$$
  in contradiction to our hypotheses. The assumption that 
   $\Phi^{*}[c_{1}(\overline{M}, \omega_{1})]= \pm 
  c_{1}(\overline{M}, \omega_{2})$ is therefore false, and the claim
  follows.
 \end{proof}
 
 Theorem 1 is now an immediate consequence, since 
 the first Chern class of a symplectic structure
is  deformation-invariant. 
 
 \begin{ack} This work was supported 
in part by  NSF grant DMS-0072591.
\end{ack}

\end{document}

%% file: gtoutput.tex
\def\ifplaintex{\expandafter\ifx\csname documentclass\endcsname\relax}

\ifplaintex 
\else
\fi

\expandafter\ifx\csname beginpicture\endcsname\relax
\expandafter\ifx\csname documentclass\endcsname\relax
\input pictex \else
\input prepictex \input pictex \input postpictex \fi\fi

\def\gt{{\mathsurround=0pt\it $\cal G\mskip-2mu$eometry \&\ 
$\cal T\!\!$opology}}        %  journal title in recommended style

\def\gtp{{\mathsurround=0pt\it $\cal G\mskip-2mu$eometry \&\ 
$\cal T\!\!$opology $\cal P\!$ublications}}  % GT publications

\def\lognumber#1{\def\thelognumber{#1}}
\def\volumenumber#1{\def\thevolumenumber{#1}}
\def\papernumber#1{\def\thepapernumber{#1}}
\def\volumeyear#1{\def\thevolumeyear{#1}}

\def\pagenumbers#1#2{\def\startpage{#1}\def\finishpage{#2}}
\def\published#1{\def\publishdate{#1}}
\def\proposed#1{\def\theproposer{#1}}
\def\seconded#1{\def\theseconders{#1}}
\def\received#1{\def\receiveddate{#1}}

\def\accepted#1{\def\accepteddate{#1}}

\long\def\asciiabstract#1{\long\def\theasciiabstract{#1}}
\def\asciikeywords#1{\def\theasciikeywords{#1}}

\let\\\par\let\thelognumber\relax
\let\thevolumenumber\relax\let\thepapernumber\relax
\let\thevolumeyear\relax\let\thesamplenumber\relax\let\startpage\relax
\let\finishpage\relax\let\publishdate\relax\let\receiveddate\relax
\let\reviseddate\relax\let\accepteddate\relax\let\theasciititle\relax
\let\theasciiauthors\relax
\let\theasciiabstract\relax\let\theasciikeywords\relax
\let\theasciiemail\relax\let\theshortauthors\relax\let\theshorttitle\relax

\long\def\maketitlep{   % start of definition of \maketitlep

\count0=\startpage

\gt\hfill      %   Journal title (top left) 
\beginpicture
\setcoordinatesystem units <0.33truein, 0.33truein> point at 2.2 0.9
\setplotsymbol ({$\cal G$})
\plotsymbolspacing=9truept
\circulararc 315 degrees from 0 1 center at 0 0
\setplotsymbol ({$\cal T$})
\circulararc 315 degrees from 1 -1 center at 1 0
\endpicture
\break
{\small\ifx\thesamplenumber\relax % sample?  
Volume \else Sample
\fi\thevolumenumber\ (\thevolumeyear)
\startpage--\finishpage\nl
Published: \publishdate}
\vglue 0.5truein plus 0.4fil minus 0.1truein

{\parskip=0pt\leftskip 0pt plus 1fil\def\\{\par\smallskip}{\ifplaintex\large
\else\Large\fi\bf\thetitle}\par\medskip}   

\vglue 0pt plus 0.1fil 

{\parskip=0pt\leftskip 0pt plus 1fil\def\\{\par}{\sc\theauthors}
\par\medskip}

\vglue 0pt plus 0.1fil 

{\small\parskip=0pt\let\newline\\
{\leftskip 0pt plus 1fil\def\\{\par}{\sl\theaddress}\par}
\expandafter\ifx\theemail\relax    % email address?
\relax\else\vglue 5pt plus 0.02fil minus 2pt\def\\{\stdspace{\rm 
and}\stdspace} 
\cl{Email:\stdspace\tt\theemail}\fi
\ifx\theurl\relax                  % URL given?
\relax\else\vglue 5pt plus 0.02fil minus 2pt\def\\{\stdspace{\rm 
and}\stdspace}
\cl{URL:\stdspace\tt\theurl}\fi\par}

\vglue 7pt plus 0.3fil minus 3pt

{\bf Abstract}
\vglue 5pt plus 0.1fil minus 2pt

\theabstract

\vglue 7pt plus 0.3fil minus 3pt

{\bf AMS Classification numbers}\quad Primary:\quad \theprimaryclass

Secondary:\quad \thesecondaryclass

\vglue 5pt plus 0.3fil minus 2pt

{\bf Keywords:}\quad \thekeywords

\vglue 10pt plus 0.5fil minus 5pt

{\small  Proposed: \theproposer\hfill Received: \receiveddate\nl
Seconded: \theseconders\hfill 
\ifx\reviseddate\relax                         % paper revised?
Accepted: \accepteddate                        % no
\else
Revised: \reviseddate                          % yes
\fi}
\eject
}       %  end of definition of \maketitlep

\let\maketitlepage\maketitlep
\let\maketitle\maketitlepage

\font\phead=cmsl9 scaled 950
\font=cmsl9 scaled 1050
\font\pnum=cmbx10 scaled 913
\font\lnum=cmbx10 
\font\pfoot=cmsl9 scaled 950
\font=cmsl9 scaled 1050
\ifplaintex
\headline{\vbox to 0pt{\vskip -4.5mm\line{\small\phead\ifnum
\count0=\startpage ISSN 1364-0380 (on line)
1465-3060 (printed) \hfill {\pnum\folio}\else\ifodd\count0\def\\{ }% 
\ifx\theshorttitle\relax\thetitle\else\theshorttitle\fi\hfill{\pnum\folio}
\else\def\\{ and }{\pnum\folio}\hfill\ifx\theshortauthors\relax\theauthors
\else\theshortauthors\fi\fi\fi}\vss}}
\footline{\vbox to 0pt{\vglue 0mm\line{\small\pfoot\ifnum\count0=\startpage
\copyright\ \gtp\hfill\else
\gt, Volume \thevolumenumber\ (\thevolumeyear)\hfill\fi}\vss
}}
\else
\makeatletter
\def\@oddhead{{\small\lhead\ifnum\count0=\startpage ISSN 1364-0380 (on line)
1465-3060 (printed) \hfill {\lnum\number\count0}\else\ifodd\count0
\def\\{ }\ifx\theshorttitle\relax \thetitle \else\theshorttitle\fi\hfill
{\lnum\number\count0}\else\def\\{ and }{\lnum\number\count0}
\hfill\ifx\theshortauthors\relax 
\theauthors\else\theshortauthors\fi\fi\fi}}\def\@evenhead{\@oddhead}
\def\@oddfoot{\small\lfoot\ifnum\count0=\startpage\copyright\ \gtp\hfill\else
\gt, Volume \thevolumenumber\ (\thevolumeyear)\hfill\fi}
\def\@evenfoot{\@oddfoot}
\makeatother
\fi

\newwrite\gtoutfile
\long\gdef\makeheadfile{  %%% start of definition of \makeheadfile
{\def\\{, }\def\s{ }
\immediate\openout\gtoutfile head.xxx
\immediate\write\gtoutfile{To: math@arxiv.org}
\immediate\write\gtoutfile{Subject: put or rep NNNNN:pppp}
\immediate\write\gtoutfile{--text follows this line--}
\immediate\write\gtoutfile{Proxy-for: \ifx\theasciiauthors\relax
\theauthors\else\theasciiauthors\fi\s<\ifx\theasciiemail\relax\theemail\else\theasciiemail\fi>}
\immediate\write\gtoutfile{\noexpand\\}
\immediate\write\gtoutfile{Authors: \ifx\theasciiauthors\relax
\theauthors\else\theasciiauthors\fi}
\immediate\write\gtoutfile{Title: \ifx\theasciititle\relax
\thetitle\else\theasciititle\fi}
\immediate\write\gtoutfile{Subj-class: GT or SG or MG etc}
\immediate\write\gtoutfile{MSC-class: \theprimaryclass\ifx\thesecondaryclass\relax\else, \thesecondaryclass\fi}
\immediate\write\gtoutfile{Journal-ref: Geom. Topol. \thevolumenumber
(\thevolumeyear) \startpage-\finishpage}
\immediate\write\gtoutfile{Comments: Published by Geometry and Topology at}
\immediate\write\gtoutfile{\s\s http://www.maths.warwick.ac.uk/gt/GTVol\thevolumenumber/paper\thepapernumber.abs.html}
\immediate\write\gtoutfile{\noexpand\\}
\immediate\write\gtoutfile{}
\ifx\theasciiabstract\relax
\immediate\write\gtoutfile{\theabstract}\else
\immediate\write\gtoutfile{\theasciiabstract}\fi
\immediate\write\gtoutfile{}
\immediate\write\gtoutfile{\noexpand\\}
\immediate\write\gtoutfile{}
\immediate\closeout\gtoutfile}}  %%% end of definition of \makeheadfile

\def\maketitlepage{\maketitlep\makeheadfile}
\let\maketitle\maketitlepage

\def\ifplaintex{\expandafter\ifx\csname documentclass\endcsname\relax}

\ifplaintex 
\else
\fi

\expandafter\ifx\csname beginpicture\endcsname\relax
\expandafter\ifx\csname documentclass\endcsname\relax
\input pictex \else
\input prepictex \input pictex \input postpictex \fi\fi

\def\gt{{\mathsurround=0pt\it $\cal G\mskip-2mu$eometry \&\ 
$\cal T\!\!$opology}}        %  journal title in recommended style

\def\gtp{{\mathsurround=0pt\it $\cal G\mskip-2mu$eometry \&\ 
$\cal T\!\!$opology $\cal P\!$ublications}}  % GT publications

\def\lognumber#1{\def\thelognumber{#1}}
\def\volumenumber#1{\def\thevolumenumber{#1}}
\def\papernumber#1{\def\thepapernumber{#1}}
\def\volumeyear#1{\def\thevolumeyear{#1}}

\def\pagenumbers#1#2{\def\startpage{#1}\def\finishpage{#2}}
\def\published#1{\def\publishdate{#1}}
\def\proposed#1{\def\theproposer{#1}}
\def\seconded#1{\def\theseconders{#1}}
\def\received#1{\def\receiveddate{#1}}

\def\accepted#1{\def\accepteddate{#1}}

\long\def\asciiabstract#1{\long\def\theasciiabstract{#1}}
\def\asciikeywords#1{\def\theasciikeywords{#1}}

\let\\\par\let\thelognumber\relax
\let\thevolumenumber\relax\let\thepapernumber\relax
\let\thevolumeyear\relax\let\thesamplenumber\relax\let\startpage\relax
\let\finishpage\relax\let\publishdate\relax\let\receiveddate\relax
\let\reviseddate\relax\let\accepteddate\relax\let\theasciititle\relax
\let\theasciiauthors\relax
\let\theasciiabstract\relax\let\theasciikeywords\relax
\let\theasciiemail\relax\let\theshortauthors\relax\let\theshorttitle\relax

\long\def\maketitlep{   % start of definition of \maketitlep

\count0=\startpage

\gt\hfill      %   Journal title (top left) 
\beginpicture
\setcoordinatesystem units <0.33truein, 0.33truein> point at 2.2 0.9
\setplotsymbol ({$\cal G$})
\plotsymbolspacing=9truept
\circulararc 315 degrees from 0 1 center at 0 0
\setplotsymbol ({$\cal T$})
\circulararc 315 degrees from 1 -1 center at 1 0
\endpicture
\break
{\small\ifx\thesamplenumber\relax % sample?  
Volume \else Sample
\fi\thevolumenumber\ (\thevolumeyear)
\startpage--\finishpage\nl
Published: \publishdate}
\vglue 0.5truein plus 0.4fil minus 0.1truein

{\parskip=0pt\leftskip 0pt plus 1fil\def\\{\par\smallskip}{\ifplaintex\large
\else\Large\fi\bf\thetitle}\par\medskip}   

\vglue 0pt plus 0.1fil 

{\parskip=0pt\leftskip 0pt plus 1fil\def\\{\par}{\sc\theauthors}
\par\medskip}

\vglue 0pt plus 0.1fil 

{\small\parskip=0pt\let\newline\\
{\leftskip 0pt plus 1fil\def\\{\par}{\sl\theaddress}\par}
\expandafter\ifx\theemail\relax    % email address?
\relax\else\vglue 5pt plus 0.02fil minus 2pt\def\\{\stdspace{\rm 
and}\stdspace} 
\cl{Email:\stdspace\tt\theemail}\fi
\ifx\theurl\relax                  % URL given?
\relax\else\vglue 5pt plus 0.02fil minus 2pt\def\\{\stdspace{\rm 
and}\stdspace}
\cl{URL:\stdspace\tt\theurl}\fi\par}

\vglue 7pt plus 0.3fil minus 3pt

{\bf Abstract}
\vglue 5pt plus 0.1fil minus 2pt

\theabstract

\vglue 7pt plus 0.3fil minus 3pt

{\bf AMS Classification numbers}\quad Primary:\quad \theprimaryclass

Secondary:\quad \thesecondaryclass

\vglue 5pt plus 0.3fil minus 2pt

{\bf Keywords:}\quad \thekeywords

\vglue 10pt plus 0.5fil minus 5pt

{\small  Proposed: \theproposer\hfill Received: \receiveddate\nl
Seconded: \theseconders\hfill 
\ifx\reviseddate\relax                         % paper revised?
Accepted: \accepteddate                        % no
\else
Revised: \reviseddate                          % yes
\fi}
\eject
}       %  end of definition of \maketitlep

\let\maketitlepage\maketitlep
\let\maketitle\maketitlepage

\font\phead=cmsl9 scaled 950
\font=cmsl9 scaled 1050
\font\pnum=cmbx10 scaled 913
\font\lnum=cmbx10 
\font\pfoot=cmsl9 scaled 950
\font=cmsl9 scaled 1050
\ifplaintex
\headline{\vbox to 0pt{\vskip -4.5mm\line{\small\phead\ifnum
\count0=\startpage ISSN 1364-0380 (on line)
1465-3060 (printed) \hfill {\pnum\folio}\else\ifodd\count0\def\\{ }% 
\ifx\theshorttitle\relax\thetitle\else\theshorttitle\fi\hfill{\pnum\folio}
\else\def\\{ and }{\pnum\folio}\hfill\ifx\theshortauthors\relax\theauthors
\else\theshortauthors\fi\fi\fi}\vss}}
\footline{\vbox to 0pt{\vglue 0mm\line{\small\pfoot\ifnum\count0=\startpage
\copyright\ \gtp\hfill\else
\gt, Volume \thevolumenumber\ (\thevolumeyear)\hfill\fi}\vss
}}
\else
\makeatletter
\def\@oddhead{{\small\lhead\ifnum\count0=\startpage ISSN 1364-0380 (on line)
1465-3060 (printed) \hfill {\lnum\number\count0}\else\ifodd\count0
\def\\{ }\ifx\theshorttitle\relax \thetitle \else\theshorttitle\fi\hfill
{\lnum\number\count0}\else\def\\{ and }{\lnum\number\count0}
\hfill\ifx\theshortauthors\relax 
\theauthors\else\theshortauthors\fi\fi\fi}}\def\@evenhead{\@oddhead}
\def\@oddfoot{\small\lfoot\ifnum\count0=\startpage\copyright\ \gtp\hfill\else
\gt, Volume \thevolumenumber\ (\thevolumeyear)\hfill\fi}
\def\@evenfoot{\@oddfoot}
\makeatother
\fi

\newwrite\gtoutfile
\long\gdef\makeheadfile{  %%% start of definition of \makeheadfile
{\def\\{, }\def\s{ }
\immediate\openout\gtoutfile head.xxx
\immediate\write\gtoutfile{To: math@arxiv.org}
\immediate\write\gtoutfile{Subject: put or rep NNNNN:pppp}
\immediate\write\gtoutfile{--text follows this line--}
\immediate\write\gtoutfile{Proxy-for: \ifx\theasciiauthors\relax
\theauthors\else\theasciiauthors\fi\s<\ifx\theasciiemail\relax\theemail\else\theasciiemail\fi>}
\immediate\write\gtoutfile{\noexpand\\}
\immediate\write\gtoutfile{Authors: \ifx\theasciiauthors\relax
\theauthors\else\theasciiauthors\fi}
\immediate\write\gtoutfile{Title: \ifx\theasciititle\relax
\thetitle\else\theasciititle\fi}
\immediate\write\gtoutfile{Subj-class: GT or SG or MG etc}
\immediate\write\gtoutfile{MSC-class: \theprimaryclass\ifx\thesecondaryclass\relax\else, \thesecondaryclass\fi}
\immediate\write\gtoutfile{Journal-ref: Geom. Topol. \thevolumenumber
(\thevolumeyear) \startpage-\finishpage}
\immediate\write\gtoutfile{Comments: Published by Geometry and Topology at}
\immediate\write\gtoutfile{\s\s http://www.maths.warwick.ac.uk/gt/GTVol\thevolumenumber/paper\thepapernumber.abs.html}
\immediate\write\gtoutfile{\noexpand\\}
\immediate\write\gtoutfile{}
\ifx\theasciiabstract\relax
\immediate\write\gtoutfile{\theabstract}\else
\immediate\write\gtoutfile{\theasciiabstract}\fi
\immediate\write\gtoutfile{}
\immediate\write\gtoutfile{\noexpand\\}
\immediate\write\gtoutfile{}
\immediate\closeout\gtoutfile}}  %%% end of definition of \makeheadfile

\def\maketitlepage{\maketitlep\makeheadfile}
\let\maketitle\maketitlepage